\documentclass[12pt,reqno]{amsart}
\usepackage{amssymb,delarray}
\usepackage{amsfonts}
\usepackage{epsfig}
\usepackage[all]{xy}
\usepackage{amscd}
\usepackage{epigraph}

\usepackage{graphicx}

\def\leq{\leqslant}
\def\geq{\geqslant}

\def\edvo{\rule {6pt}{6pt}} 

\newtheorem{thm}{Theorem}
\newtheorem{lem}{Lemma}
\newtheorem{cor}{Corollary}
{\catcode`\@=11
\gdef\n@te#1#2{\leavevmode\vadjust{%
 {\setbox\z@\hbox to\z@{\strut#1}%
  \setbox\z@\hbox{\raise\dp\strutbox\box\z@}\ht\z@=\z@\dp\z@=\z@%
  #2\box\z@}}}
\gdef\leftnote#1{\n@te{\hss#1\quad}{}}
\gdef\rightnote#1{\n@te{\quad\kern-\leftskip#1\hss}{\moveright\hsize}}
\gdef\?{\FN@\qumark}
\gdef\qumark{\ifx\next"\DN@"##1"{\leftnote{\rm##1}}\else
 \DN@{\leftnote{\rm??}}\fi{\rm??}\next@}}

\begin{document}

\author{O.V.~Kulikova}
\email{olga.kulikova@mail.ru}

\title{On independent families  of normal subgroups in free groups}

\markboth{}{On independent families of normal subgroups in free groups}

\maketitle

\begin{abstract}
Consider a presentation $\mathcal{P}=<{\bf x}\mid{\bf \bigcup_{i=1}^n r_i}>$.
 Let ${\bf R_i}$ be the normal closure of the set ${\bf r_i}$ in the free group ${\bf F}$ with basis ${\bf x}$, $\mathcal{P}_i=<{\bf x}\mid{\bf r_i}>$, ${\bf N_i} = \prod_{j\neq i}{\bf R_j}$. In the present article, using geometric techniques of pictures, generators for  $\frac{{\bf R_i}\cap {\bf N_i}}{[{\bf R_i}, {\bf N_i}]}$, $i=1,\ldots,n$, are obtained from a set of generators over $\{\mathcal{P}_i\mid i=1,\ldots, n\}$ for  $\pi_2(\mathcal{P})$. As a corollary, we get a sufficient condition for the family $\{{\bf R_1},\ldots,{\bf R_n}\}$ to be independent.
\end{abstract}

\section*{Introduction}
Consider a presentation $\mathcal{P}=<{\bf x}\mid{\bf r}>$, where ${\bf r=\bigcup_{i=1}^n r_i}$.
For $i=1,\ldots, n$, let ${\bf R_i}$ be the normal closure of the set ${\bf r_i}$ in the free group ${\bf F}$ with basis ${\bf x}$, $\mathcal{P}_i=<{\bf x}\mid{\bf r_i}>$, ${\bf N_i} = \prod_{j\neq i}{\bf R_j}$.
Set ${\bf R} = \prod_{i=1}^n{\bf R_i}$, ${\bf G} = {\bf F}/{\bf R}$  and ${\bf G_i} = {\bf F}/{\bf R_i}$.
In this paper we study the quotient groups $\frac{{\bf R_i}\cap {\bf N_i}}{[{\bf R_i}, {\bf N_i}]}$, $i=1,\ldots,n$, which are a natural measure of the redundancy between ${\bf R_1}, ..., {\bf R_n}$.

Geometric techniques of spherical pictures \cite{pride, gener_b_pr} are used to prove the main theorem of this paper, in which we determine generators for $\frac{{\bf R_i}\cap {\bf N_i}}{[{\bf R_i}, {\bf N_i}]}$, $i=1,\ldots,n$, providing that a set of spherical pictures over $\mathcal{P}$, generating $\pi_2(\mathcal{P})$ over $\{\mathcal{P}_i\mid i=1,\ldots, n\}$, is known.  The idea of applying these techniques to determine generators was already used in \cite{BoGu92} for some presentations.

The family $\{{\bf R_1},\ldots,{\bf R_n}\}$ is said to be {\it independent}, if ${\bf R_i}\cap {\bf N_i} = [{\bf R_i}, {\bf N_i}]$ for $i=1,\ldots,n$. Independence may be considered as ensuring that certain intersections are as small as possible, or as ensuring that certain commutator subgroups are as large as possible. Independence and related notions have been studied in \cite{Bo91, ratcl, DuElGi92, gilbert, huebsch, lind=my}. For example, it is shown in \cite{DuElGi92} that $\{{\bf R_1},\ldots,{\bf R_n}\}$ is independent if and only if the inclusions ${\bf R_i}\rightarrow{\bf R}$ induce an isomorphism $\oplus_{i=1}^n(\mathbb{Z}{\bf G}\otimes_{\mathbb{Z}{\bf G_i}}H_1{\bf R_i})\rightarrow H_1{\bf R}$ of relation modules. Also it is known \cite{ratcl, DuElGi92, huebsch} that if $\{{\bf R_1},\ldots,{\bf R_n}\}$ is independent, then $\mathcal{P}$ is (A)  over $\{\mathcal{P}_i\mid i=1,\ldots,n\}$, i.e. $\pi_2(\mathcal{P})$ is generated by the empty set over $\{\mathcal{P}_i\mid i\in I\}$. It is proved in \cite{ratcl} that the converse holds in case $n=2$. In the present paper we obtain the converse statement for arbitrary $n$, as a corollary of the main theorem.

In the first section of this paper we remind  the facts from  \cite{pride} and \cite{gener_b_pr}, required to formulate and prove the main results of this paper. The main results themselves are formulated in the second section and  proved in the third section. The reader may omit the first section and skip to the next ones if the papers \cite{pride} and \cite{gener_b_pr} are familiar to him.

\section{Main definitions}\label{S:df_in_F}
The following definitions, notation and facts from \cite{pride} and \cite{gener_b_pr} will be used in the second and third sections of the present paper.

Let $\mathcal{P} = \langle \bf{x}\mid \bf{r}\rangle$ be a presentation for a group $\bf G$, where $\bf{x}$  is a set ("generators") and $ \bf{r}$ is a set of cyclically reduced words on $\bf{x}\cup\bf{x}^{-1}$ ("relators").
We assume that no relator of $\bf{r}$ is freely trivial, nor is a conjugate of any other relator or its inverse (so called  the $RH$-hypothesis).

We let $\bf w$ denote the set of all words on $\bf{x}\cup\bf{x}^{-1}$ (reduced or not). Two words $u$ and $v$ in $\bf w$ are {\it identically equal}, if they
represent the same element in the free semigroup on $\bf{x}\cup\bf{x}^{-1}$. The words  $u$ and $v$ are {\it freely equal}, if one can be obtained from the other by a finite number of insertions and deletions of inverse pairs $x^{\varepsilon}x^{-\varepsilon}$, where $x\in{\bf x}, \varepsilon=\pm 1$.
The free equivalence class of $W\in\bf{w}$ will be denoted by $[W]$.
The free group  ${\bf F}$  on $\bf{x}$ consists of the free equivalence classes, where the multiplication is defined by $[W][U]=[WU]$. We let ${\bf N}$ denote  the normal closure of $\{[R]\mid R\in\bf{r}\}$ in ${\bf F}$.
Thus, ${\bf G} = {\bf F}/ {\bf N}$.

If $\bf s$ is a subset of $\bf r$ then we let $\bf s^w$ denote the set of all words of the form $WS^{\epsilon}W^{-1}$, where $W\in{\bf w}, S\in{\bf s}, \epsilon=\pm 1$.

\underline{Sequences.}

Let us consider finite sequences of elements of $\bf r^w$.

Let $\sigma = (c_1, ...., c_m)$, where $c_i\in \bf r^w$ ($i=1,...,m$). We define the {\it inverse} $\sigma^{-1}$ of $\sigma$ to be $(c_m^{-1},...,c_1^{-1})$. We define the {\it conjugate} $W\sigma W^{-1}$ of $\sigma$ by $W\in {\bf w}$ to be $(Wc_1W^{-1},...,Wc_mW^{-1})$.

We define operations on sequences as follows. Let $\sigma = (c_1, ...., c_m)$, where $c_i = W_iR_i^{\epsilon_i}W_i^{-1}$,  $W_i\in{\bf w}, R_i\in{\bf r}, \epsilon_i=\pm 1, i=1,...,m$.

(SUB) {\it Substitution.} Replace each $W_i$ by a word freely equal to it.

(DEL) {\it Deletion.} Delete two consecutive terms if one is identically equal to the inverse of the other.

(INS) {\it Insertion. } The opposite of deletion.

(EX) {\it Exchange. } Replace two consecutive terms $c_i$, $c_{i+1}$ by either $c_{i+1}$, $c_{i+1}^{-1}c_ic_{i+1}$ or by $c_ic_{i+1}c_i^{-1}$, $c_i$.

Two sequences $\sigma$, $\sigma'$ will be said to be {\it (Peiffer) equivalent} (denoted $\sigma\sim\sigma'$) if one can be obtained from the other by a finite number of applications of the operations (SUB), (DEL), (INS), (EX). The equivalence class containing $\sigma$ will be denoted by $<\sigma>$.

We can define a binary operation $+$ on the set $\Sigma$ of all equivalence classes by the rule
$$
<\sigma_1> + <\sigma_2> = <\sigma_1\sigma_2>.
$$
(Here $\sigma_1\sigma_2$ is the juxtaposition of the two sequences $\sigma_1$, $\sigma_2$.) Under this operation $\Sigma$ is a group. The identity (zero element) is the equivalence class of the empty sequence, and the inverse (negative) $-<\sigma>$ of $<\sigma>$ is $<\sigma^{-1}>$.

We define $\Pi\sigma$ to be the product $c_1c_2...c_m$. We say that $\sigma$
is an {\it identity sequence}, if  $\Pi\sigma$ is freely equal to $1$. We let $\pi_2$ denote the (abelian) subgroup of $\Sigma$ consisting of all elements  $<\sigma>$, where $\sigma$ is an identity sequence. (Occasionally, when we want to emphasize the presentation $\mathcal{P}$, we will write $\pi_2(\mathcal{P})$).

${\bf F}$ acts on $\Sigma$ by the rule
$[W]\cdot<\sigma> = <W\sigma W^{-1}>$,  where $[W]\in {\bf F}, <\sigma>\in \Sigma$.
The mapping
$\partial : \Sigma\rightarrow {\bf F}$, where
$\partial : <\sigma>\mapsto[\Pi\sigma]$,
is a group homomorphism  with the image ${\bf N}$ and the kernel $\pi_2$.
The triple $(\Sigma, {\bf F}, \partial)$ is an example of  a (free) crossed module.
${\bf N}$ acts trivially on $\pi_2$. It follows that $\pi_2$ is a  left ${\bf G}$-module with ${\bf G}$-action given by
$[W]{\bf N}\cdot<\sigma> = [W]\cdot<\sigma>$, where $[W]\in {\bf F}$, $<\sigma>\in \pi_2$.

 Let $R\in {\bf r}$. Then $R= {R^o}^{p(R)}$, where $R^o$ is not a proper power and $p(R)$ is a positive integer ($R^o$ is the {\it root} of $R$, and $p(R)$ is the {\it period}).
The identity sequences
$$
\zeta_R = (R, R^oR^{-1}{R^o}^{-1})\,\,\,\,\,\, (R\in {\bf r})
$$
will be called the {\it trivial sequences}, and the elements
$$
<\zeta_R>
$$
of $\pi_2$ will be called the {\it trivial elements}. The submodule of $\pi_2$ generated by the trivial elements will be denoted by ${\bf T}$.

$\mathcal{P}$ is called {\it aspherical} (respectively, {\it combinatorially aspherical (CA)}) if $\pi_2=0$ (respectively, $\pi_2={\bf T}$).
As is shown in \cite{cch}, $\mathcal{P}$ is aspherical if and  only if $\mathcal{P}$ is (CA) and no element of $\bf r$ is a proper power.

\underline{Subpresentations.}

Let $\mathcal{P}_0 = \langle \bf{x}_0\mid \bf{r}_0\rangle$ be a subpresentation of $\mathcal{P}$, and let  $(\Sigma_0, {\bf F}_0, \partial_0)$ be the associated crossed module. There is an obvious mapping of crossed modules
$$
(\varphi,\psi) : (\Sigma_0, {\bf F}_0, \partial_0)\rightarrow (\Sigma, {\bf F}, \partial),
$$
where
$$
\varphi(<\sigma>_0) = <\sigma>\,\,\,\,\,\, (<\sigma>_0\in \Sigma_0);
$$
$$
\psi([W]_0)=[W]\,\,\,\,\,\, ([W]_0\in {\bf F}_0).
$$
Restricting $\varphi$ gives a homomorphism
$$
j: \pi_2(\mathcal{P}_0)\rightarrow\pi_2(\mathcal{P}).
$$
In general $j$ is not injective.
It is still unknown whether
$j$ is injective for every subpresentation  $\mathcal{P}_0$ of aspherical $\mathcal{P}$ (Whitehead problem).

Consider the presentation  $\mathcal{P} = \langle \bf{x}\mid \bf{r}\rangle$ and suppose  that $\bf r$ is expressed as a union
$\bf r = r_1\cup r_2$. For $\lambda =1,2$, let  $\mathcal{P}_{\lambda} = \langle \bf{x}\mid \bf{r}_{\lambda}\rangle$ and $j_{\lambda}: \pi_2(\mathcal{P_{\lambda}})\rightarrow\pi_2(\mathcal{P})$ be the homomorphism as discussed above. Note that $\rm{Im}\, j_{\lambda}$ is a submodule of
$\pi_2(\mathcal{P})$.

Let ${\bf N_{\lambda}}$, where $\lambda = 1,2$, be the normal closure of $\{[R]\mid R\in {\bf r}_{\lambda}\}$ in the free group ${\bf F}$. Now ${\bf F}$ acts on $\frac{{\bf N_1}\cap {\bf N_2}}{[{\bf N_1},{\bf N_2}]}$ via conjugation:
$$
[W]\cdot[U][{\bf N_1}, {\bf N_2}] = [WUW^{-1}][{\bf N_1}, {\bf N_2}]\,\,\,\,\,\, ([W]\in {\bf F}, [U]\in {\bf N_1}\cap {\bf N_2}).
$$
It is easy to show that ${\bf N}(={\bf N_1}{\bf N_2})$ acts trivially, and so we get an induced action of ${\bf G} = {\bf F}/{\bf N}$  on $\frac{{\bf N_1}\cap {\bf N_2}}{[{\bf N_1},{\bf N_2}]}$.
We can define a ${\bf G}$-homomorphism
$$
\eta: \pi_2(\mathcal{P})\rightarrow \frac{{\bf N_1}\cap {\bf N_2}}{[{\bf N_1},{\bf N_2}]}
$$
by the following rule. Let $<\sigma>\in \pi_2(\mathcal{P})$, and $V$ be the product (taken in order) of the elements of $\sigma$ which belong to $\bf r_1^w$. Then $\eta(<\sigma>) = [V][{\bf N_1}, {\bf N_2}]$. It is not hard to  show that $\eta$ is well-defined.

The following theorem was proved by Gutierr$\acute{\text е}$z and Ratcliffe \cite{ratcl}.

\begin{thm} \label{thmGR}
 Let $\xi: \rm{Im}\, j_1\oplus\rm{Im}\, j_2 \longrightarrow \pi_2(\mathcal{P})$ be induced by the inclusions $\rm{Im}\, j_1,\rm{Im}\, j_2 \longrightarrow \pi_2(\mathcal{P})$. Then the sequence
$$
\rm{Im}\, j_1\oplus\rm{Im}\, j_2 \longrightarrow^{\xi} \pi_2(\mathcal{P})\longrightarrow^{\eta} \frac{{\bf N_1}\cap {\bf N_2}}{[{\bf N_1},{\bf N_2}]}\longrightarrow 0
$$
is exact. If $\bf r_1$ and $\bf r_2$ are disjoint then $\xi$ is injective, and so the sequence is short exact.
\end{thm}

\underline{Pictures.}

Sequences can be studied very successfully using geometric objects called pictures. Pictures were introduced in \cite{igusa, rourke}. These objects are  a very useful tool to solve combinatorial group theory problems. See, for example,  \cite  {pride, gener_b_pr} and the references cited there.

Let us remind the definition of pictures (according to \cite{gener_b_pr}).

A {\it picture} $\mathbf{P}$ is a finite collections of pairwise disjoint closed disks $\Delta_1,\ldots, \Delta_m$ in a closed disk $D^2$,
together with a finite number of disjoint arcs $\alpha_1,\ldots,\alpha_l$ properly embedded in $D^2 -\cup_{i=1}^mint\Delta_i$
(where "$int$"\, denotes interior). Loosely speaking, below the disks $\Delta_1,\ldots,\Delta_m$ will be called {\it vertices} of $\mathbf{P}$. An arc
can be either a simple closed curve having trivial intersection  with
с $\partial D^2\cup\partial\Delta_1\cup\ldots\cup\partial\Delta_m$, or a simple non-closed curve which joins two different points of $\partial D^2\cup\partial\Delta_1\cup\ldots\cup\partial\Delta_m$.
The {\it boundary} $\partial D^2$ of $\mathbf{P}$ will be denoted by $\partial \mathbf{P}$.
The {\it corners} of a vertex $\Delta$ of $\mathbf{P}$ are the closures of the connected components of $\partial\Delta -\cup_{j}\alpha_j$.
The {\it regions} of $\mathbf{P}$ are the closures of the connected components of $D^2-(\cup_i\Delta_i\bigcup\cup_j\alpha_j)$.
The {\it components} of $\mathbf{P}$ are the connected components of $\cup_i\Delta_i\cup\cup_j\alpha_j$.
The picture $\mathbf{P}$ is {\it connected} if it has at most one component. The picture $\mathbf{P}$ is {\it spherical} if it has
at least one vertex and no arc of $\mathbf{P}$ meets $\partial \mathbf{P}$.

Assume that a group presentation $\mathcal{P} = \langle \bf{x}\mid \bf{r}\rangle$ and a picture $\mathbf{P}$ are given. Fix an
orientation of the ambient disk $D^2$, thereby determining a sense of positive rotation (i.g. clockwise). Assume that the vertices and arcs of $\mathbf{P}$ are labeled by elements of $\mathcal{P}$ as follows.

(i) Each arc of $\mathbf{P}$ is equipped with a normal orientation (indicated  by an arrow transverse to the arc), and is labeled by an element of $\bf{x}\cup\bf{x}^{-1}$.

(ii) Each vertex $\Delta$ of $\mathbf{P}$ is equipped with a sign $\epsilon(\Delta) = \pm 1$ and is labeled by a relator $R(\Delta)\in\bf{r}$.

For a corner $c$ of a vertex $\Delta$ of $\mathbf{P}$, $W(c)$ denotes the word in $\bf{x}\cup\bf{x}^{-1}$ obtained by reading in order the (signed) labels on the arcs that are encountered in a walk around $\partial\Delta$ in the positive direction, beginning and ending at an interior point of $c$ (with the understanding that if we cross an arc, labeled $y$ say, in the direction of its normal orientation then we read $y$, whereas  if we cross the arc against the orientation we read $y^{-1}$).
The oriented and labeled picture $\mathbf{P}$ is a {\it picture over} $\mathcal{P}$ if for each corner $c$ of each vertex $\Delta$ of $\mathbf{P}$, $W(c)$ is identically equal to a  cyclic permutation of $R(\Delta)^{\epsilon(\Delta)}$. We call  $W(c)$ the {\it label} of $\Delta$, and denote it by $W(\Delta)$.

A corner $c$ is a {\it basic corner} of $\Delta$ of $\mathbf{P}$ if $W(c)$ is identically equal to $R(\Delta)^{\epsilon(\Delta)}$.  The vertex $\Delta$ has exactly $p$ basic corners, where $p\geq 1$ is the period of $R(\Delta)$.

A picture $\mathbf{P}$ over $\mathcal{P}$ becomes a {\it based} picture over $\mathcal{P}$ when it is equipped with basepoints as follows.

\begin{itemize}
\item Each vertex $\Delta$ has one {\it basepoint}, which is a selected point in the interior of a basic corner of $\Delta$.

\item $\mathbf{P}$ has a {\it global basepoint}, which is a selected point in $\partial \mathbf{P}$ that does not lie on any arc of $\mathbf{P}$.

\end{itemize}

The {\it boundary label} on a based picture $\mathbf{P}$ over $\mathcal{P}$ is the word $W(\mathbf{P})$ obtained by reading in order the (signed) labels on the arcs of $\mathbf{P}$ that are encountered in a walk around $\partial \mathbf{P}$ in the positive direction, beginning and ending at the global basepoint. Alteration of the global basepoint or of the orientation of the ambient disk $D^2$ changes the boundary label of $\mathbf{P}$ only up to cyclic permutation and inversion.

There is the following pictorial version of the "van Kampen lemma"\, (it can be obtained from
the theorem 1.1 (V) and the lemma 1.2 (V) of \cite{lind} and duality).

\begin{lem}\label{lemVKampen} A word $U$ in $\bf{x}\cup\bf{x}^{-1}$ represents the identity of ${\bf G}$ defined by $\mathcal{P} = \langle \bf{x}\mid \bf{r}\rangle$ if and only if there is a based picture $\mathbf{P}$ over $\mathcal{P}$ with boundary label identically equal to $U$.
\end{lem}

The mirror image of a picture $\mathbf{P}$ will be denoted by $-\mathbf{P}$.
We can form the {\it sum} $\mathbf{P}_1+\mathbf{P}_2$ of two pictures $\mathbf{P}_1$ and $\mathbf{P}_2$ in the obvious way (Figure 1):

\unitlength=1mm
\begin{picture}(120,50)(-20,0)

       \put(20,25){\circle{20}}
       \put(20,18){\circle*{1}}
       \put(19,24){$\mathbf{P}_1$}
       \put(40,25){+}
       \put(60,25){\circle{20}}
       \put(60,18){\circle*{1}}
       \put(59,24){$\mathbf{P}_2$}
       \put(80,25){=}
       \put(100,25){\circle{20}}
       \put(100,18){\circle*{1}}
       \multiput(100,19)(0,2){7}{\line(0,1){1}}
       \put(95,24){$\mathbf{P}_1$}
       \put(101,24){$\mathbf{P}_2$}
\end{picture}

\begin{center}
Figure 1
\end{center}

A {\it transverse path $\gamma$} in $\mathbf{P}$ over $\mathcal{P}$ is a path in the closure of $D^2-\bigcup_i \Delta_i$ which intersects the arcs of $\mathbf{P}$ only finitely many times (moreover, if the path intersects an arc then it crossed it, and doesn't just touch it), no endpoints of $\gamma$ touches any arc, and whenever $\gamma$ meets  $\partial\mathcal{P}$ or any $\partial\Delta_i$, it does so only in its endpoints. Since we will only over consider transverse paths, we will from now on drop the use of the adjective "transverse".

If we travel along a path $\gamma$ from its initial point to its terminal point  we will cross various arcs, and we can read off the (signed) labels on these arcs, giving a word $W(\gamma)$, the {\it label on} $\gamma$.

Let a simple closed path $\gamma$ in $\mathbf{P}$ encloses a {\it subpicture} $\mathbf{B}$ of $\mathbf{P}$. This subpicture consists of the vertices and (portions of) arcs that are separated from $\partial\mathbf{P}$ by $\gamma$. When $\mathbf{P}$ is spherical, the {\it compliment} of $\mathbf{B}$
in $\mathbf{P}$ is defined as follows. Delete the interior of $\mathbf{B}$ to form an oriented annulus. Identification of $\partial\mathbf{P}$ to a point produces an oriented disk that has the boundary $\gamma$, and which supports a new picture over $\mathcal{P}$. The compliment of $\mathbf{B}$
 in $\mathbf{P}$ is obtained from this new picture by a planar reflection. The complement has the same boundary label as $\mathbf{B}$ and its vertices are those of $\mathbf{P}-\mathbf{B}$, taken with opposite signs.

 The subpicture $\mathbf{B}$ enclosed by a simple closed path $\gamma$ will be called a {\it spherical subpicture} if $\gamma$ intersects no arc. A spherical subpicture will be called {\it empty} if it neither consists of any vertex nor any portion of any arc.

\underline{On the connection between sequences and pictures.}

A {\it spray} for a based picture $\mathbf{P}$ with $n$ vertices $\Delta_1,\Delta_2,\ldots,\Delta_n$ is a sequence $\underline{\gamma}=(\gamma_1, \gamma_2, \ldots, \gamma_n)$ of simple paths satisfying the following:

\begin{itemize}
\item for each $i=1, 2,\ldots,n$, $\gamma_i$ starts at the global basepoint of $\mathbf{P}$ and ends at a basepoint of $\Delta_{\vartheta(i)}$, where $\vartheta$ is a permutation of $\{1, 2, \ldots,n\}$ (depending on $\underline{\gamma}$);

\item for $1\leq i<j\leq n$, $\gamma_i$ and $\gamma_j$ intersect only at the global basepoint;

\item travelling around the global basepoint clockwise we encounter the paths in the order $\gamma_1, \gamma_2, \ldots, \gamma_n$.
\end{itemize}
The {\it sequence $\sigma(\underline{\gamma})$ associated with $\underline{\gamma}$} is
$$(W(\gamma_1)W(\Delta_{\vartheta(1)})W(\gamma_1)^{-1}, \ldots, W(\gamma_n)W(\Delta_{\vartheta(n)})W(\gamma_n)^{-1}).$$
A based picture will be said to {\it represent} a sequence $\sigma$ if there is a spray for the picture whose associated sequence is $\sigma$. Note that if $\mathbf{P}$ represents $\sigma$ then $-\mathbf{P}$ represents $\sigma^{-1}$; if $\mathbf{P}_1$, $\mathbf{P}_2$ represents $\sigma_1$, $\sigma_2$ respectively then $\mathbf{P}_1 + \mathbf{P}_2$ represents $\sigma_1\sigma_2$. One can prove (see for example \cite{pride}) that every sequence can be represented by a picture, and every identity sequence can be represented by a spherical picture; conversely, if $\mathbf{P}$ is a picture and if $\underline{\gamma}$ is a spray for $\mathbf{P}$, then $$\partial(<\sigma(\underline{\gamma})>) = [W(\mathbf{P})],$$ and if $\mathbf{P}$ is a spherical picture and and if $\underline{\gamma}$ is a spray for $\mathbf{P}$, then $\sigma(\underline{\gamma})$ is an identity sequence. If $\underline{\gamma}$, $\underline{\gamma}'$ are two sprays for a picture $\mathbf{P}$, then $<\sigma(\underline{\gamma})> - <\sigma(\underline{\gamma}')>\in {\bf T}$ (theorem 1.4, theorem 2.4 of \cite{pride}).

Consider a set ${\bf X} = \{\mathbf{P}_{\lambda}\mid \lambda\in \Lambda\}$ of based spherical pictures over $\mathcal{P}$.
For each $\lambda$, let $\sigma_{\lambda}$ be the identity sequence associated with a spray for $\mathbf{P}_{\lambda}$, and $J({\bf X})$ be the submodule of $\pi_2$ generated by $\{<\sigma_{\lambda}>\mid \lambda\in \Lambda\}$. We say that {\it ${\bf X}$ generates
$\pi_2$} if $\pi_2 = J({\bf X}) + {\bf T}$.

It follows from Theorem \ref{thmGR} that if $\pi_2 = {\bf T}$, i.e. the presentation is (CA), then $\frac{{\bf N_1}\cap {\bf N_2}}{[{\bf N_1}, {\bf N_2}]}=0$, since $\eta({\bf T}) = 0$.
If $\pi_2 = J({\bf X}) + {\bf T}$ then $\frac{{\bf N_1}\cap {\bf N_2}}{[{\bf N_1}, {\bf N_2}]}$ normally generated by the images of the elements $<\sigma_{\lambda}>$ associated with sprays for all pictures $\mathbf{P}_{\lambda}\in {\bf X}$. Moreover, as noted in \cite{BoGu92}, the image of $<\sigma_{\lambda}>$ under $\eta$ is $[V_{\lambda}][{\bf N_1}, {\bf N_2}]$, where $V_{\lambda}$ is the label of a simple closed path in $\mathbf{P}_{\lambda}$ (oriented as $\partial\mathbf{P}_{\lambda}$) separating the vertices  with $\bf r_1$-labels from the vertices with $\bf r_2$-labels.

\underline{Operations on pictures.}

Generally below we will not distinguish between pictures which are isotopic.

The following operations ("deformations") can be applied to a based picture $\mathbf{P}$ over $\mathcal{P}$ (\cite{gener_b_pr}).

$BRIDGE:$ (Bridge move) See Figure 2.

\unitlength=1mm
\begin{picture}(120,50)(-8,0)
       \thicklines

       \put(70,25){\vector(1,0){10}}
       \thinlines
       \put(40,40){\line(0,-1){30}}
       \put(50,10){\line(0,1){30}}
       \put(38,25){\vector(1,0){5}}
       \put(52,25){\vector(-1,0){5}}

       \multiput(95,17)(0,23){2}{\line(0,-1){7}}
       \multiput(105,10)(0,23){2}{\line(0,1){7}}
       \put(100,33){\oval(10,10)[b]}
       \put(100,17){\oval(10,10)[t]}
       \put(100,24){\vector(0,-1){5}}
       \put(100,26){\vector(0,1){5}}

       \multiput(35,25)(18,0){2}{\rm $x$}
       \put(96,22){\rm $x$}
       \put(102,26){\rm $x$}
\end{picture}

\begin{center}
Figure 2
\end{center}

$FLOAT:$ Deletion of a closed arc  that separates $D^2$ into two parts, one of which contains the global basepoint of $\mathbf{P}$ and all remaining arcs and vertices of  $\mathbf{P}$ (such a closed arc is called a  {\it floating circle}).

$FLOAT^{-1}:$ (Insertion of a floating circle). The opposite of  $FLOAT$.

A {\it folding pair} is a connected spherical subpicture that contains exactly two vertices such that
\begin{itemize}
\item these two vertices are labeled by the same relator and have opposite signs;

\item the basepoints of the vertices lie in the same region;

\item each arc in the subpicture has an endpoint on each vertex.

\end{itemize}

$FOLD:$ (Deletion of a folding pair). If there is a simple closed path $\beta$ in $D^2$ such that the part of $\mathbf{P}$ encircled by $\beta$ is a folding pair, then remove that part of $\mathbf{P}$ encircled by $\beta$.

$FOLD^{-1}:$ (Insertion of a folding pair). The opposite of $FOLD$.

Let ${\bf X} = \{\mathbf{P}_{\lambda}\mid \lambda\in \Lambda\}$ be a set of based spherical pictures  over $\mathcal{P}$. By an {\it $\bf X$-picture} we mean either a picture $\mathbf{P}_{\lambda}$ from ${\bf X}$ or its mirror image $-\mathbf{P}_{\lambda}$.

$REPLACE({\bf X}):$ Replace a subpicture of $\mathbf{P}$ by the complement of that subpicture in an ${\bf X}$-picture.

Two based spherical pictures are called {\it ${\bf X}$-equivalent} if one of them can be transformed into the other one (up to planar isotopy) by a finite sequence of operations $BRIDGE$, $FLOAT^{\pm 1}$, $FOLD^{\pm 1}$, $REPLACE({\bf X})$.

We introduce two further operations on pictures as follows.

$DELETE({\bf X}):$ (Deletion of an $\bf X$-picture). If there is a simple closed path $\beta$ in $D^2$ such that the part of $\mathbf{P}$ enclosed by $\beta$ is a copy of an  ${\bf X}$-picture, then delete that part of $\mathbf{P}$ enclosed by $\beta$.

$DELETE({\bf X})^{-1}:$ (Insertion of an $\bf X$-picture). The opposite of $DELETE({\bf X})$.

Note that the operation $REPLACE({\bf X})$ includes $DELETE({\bf X})^{\pm 1}$.
On the other hand, the result of the operation $REPLACE({\bf X})$ can be obtained by a finite sequence of operations $DELETE({\bf X})^{\pm 1}$, $BRIDGE$, $FLOAT^{\pm 1}$, $FOLD^{\pm 1}$. Thus, in the definition of ${\bf X}$-equivalent pictures, $REPLACE({\bf X})$ can be replaced by $DELETE({\bf X})^{\pm 1}$.

A {\it dipole} in a picture over $\mathcal{P}$ consists of an arc which meets two corners $c_1$, $c_2$ in distinct vertices such that
\begin{itemize}
\item the two vertices are labeled by the same relator and have opposite signs;
\item $c_1$ and $c_2$ lie in the same region of the picture;
\item $W(c_1) = W(c_2)^{-1}$.
\end{itemize}

By a {\it complete dipole} over $\mathcal{P}$, we mean a connected based spherical picture over $\mathcal{P}$ that contains just two vertices, and where each arc of the picture constitutes a dipole. Note that a complete dipole is just a folding pair, in that  the vertex basepoints need not lie in the same region. If the relator that labels the two vertices  of the complete dipole has period one, then a complete dipole is exactly the same as a folding pair.  A complete dipole will be called {\it primite} if the relator labeling its vertices  has root $Q$ and period $p>1$, and there is a path joining the vertex basepoints with label $Q^f$, where {\rm gcd} $(f,p)=1$.
It follows from Lemma 2.1 \cite{gener_b_pr} that, modulo primitive dipoles, one need not be concerned with choices of vertex basepoints.

The following theorem (see Corollary 1 of Theorem 2.6 \cite{pride}, Theorem 1.6 (2) \cite{gener_b_pr}) will play an important role in the proof of the main theorem of the present paper.

\begin{thm} \label{thmPr3} Let $\bf X_0$ be a collection of based spherical pictures. Then ${\bf X_0}$ generates $\pi_2 (\mathcal{P})$
if and only if every spherical picture over  $\mathcal{P}$ is $\bf X$-equivalent to the empty picture, where $\bf X$ is the union of $\bf X_0$ and the collection of primitive dipoles for all relators of $\mathcal{P}$, which are a proper power.
\end{thm}

It follows from Theorem \ref{thmPr3} (see Corollary 2 of Theorem 2.6 \cite{pride}) that
 $\mathcal{P}$ is (CA) (i.e. $\pi_2 = {\bf T}$) if and only if every spherical picture over $\mathcal{P}$ is $\bf X$-equivalent to the empty picture, where $\bf X$ is the collection of primitive dipoles for all relators of $\mathcal{P}$, which are a proper power.

\underline{Independent sets.}

Suppose that $\{\mathcal{P}_i\mid i\in I\}$ is a collection of subpresentations of $\mathcal{P}$, and let ${\bf X}_i$ denote the collection of all based spherical pictures over $\mathcal{P}_i$, $i\in I$. We shall say that a set ${\bf Y}$ of based spherical pictures over $\mathcal{P}$ {\it generates} $\pi_2(\mathcal{P})$ {\it over } $\{\mathcal{P}_i\mid i\in I\}$ if the ${\bf G}$-module $\pi_2(\mathcal{P})$ is generated by the homotopy classes $[f_\mathbf{P}]$ ($\mathbf{P}\in{\bf Y}\cup\bigcup_{i\in I}{\bf X}_i$) (\cite{gener_b_pr}).
By the analogue of Theorem \ref{thmPr3} (see Theorem 1.6 (2) \cite{gener_b_pr}), a collection ${\bf Y}$ of based spherical pictures over $\mathcal{P}$ generates  $\pi_2(\mathcal{P})$ over  $\{\mathcal{P}_i\mid i\in I\}$ if and only if every spherical picture over  $\mathcal{P}$ is $\bf X$-equivalent to the empty picture, where  ${\bf X} = {\bf Y}\cup\bigcup_{i\in I}{\bf X}_i$.

 This notion is useful when there are certain subpresentations in the presentation which we know little about, or which in some way are arbitrary. Then, we can often isolate them  by determining a "nice"\, set of generators of $\pi_2(\mathcal{P})$ relative to these subpresentations.  Examples of sets of generators over subpresentations for $\pi_2(\mathcal{P})$ can be found in \cite{gener_b_pr}.

We will say that $\mathcal{P}$ is (CA) {\it over}  $\{\mathcal{P}_i\mid i\in I\}$ if  $\pi_2(\mathcal{P})$ is generated over $\{\mathcal{P}_i\mid i\in I\}$ by primitive dipoles; $\mathcal{P}$ is (A) {\it over}  $\{\mathcal{P}_i\mid i\in I\}$ if  $\pi_2(\mathcal{P})$ is generated over $\{\mathcal{P}_i\mid i\in I\}$ by the empty set.

Consider a presentation $\mathcal{P}=<{\bf x}\mid{\bf r}>$, where ${\bf r=\bigcup_{i=1}^n r_i}$.
For $i=1,\ldots, n$, let ${\bf R_i}$ be the normal closure of ${\bf r_i}$ in the free group ${\bf F}$ with basis ${\bf x}$, $\mathcal{P}_i=<{\bf x}\mid{\bf r_i}>$, ${\bf N_i} = \prod_{j\neq i}{\bf R_j}$.
Set ${\bf R} = \prod_{i=1}^n{\bf R_i}$, ${\bf G} = {\bf F}/{\bf R}$  и ${\bf G_i} = {\bf F}/{\bf R_i}$. The family $\{ {\bf R_1},\ldots,{\bf R_n}\}$ is said to be {\it independent} if ${\bf R_i}\cap {\bf N_i} = [{\bf R_i}, {\bf N_i}]$ for $i=1,\ldots,n$.
This and related notions have been studied in \cite{Bo91, ratcl, DuElGi92, gilbert,  huebsch}.
Note that any primitive dipole over $\mathcal{P}$ belongs to some collection ${\bf X}_i$ of all based spherical pictures  over $\mathcal{P}_i$ $(i=1,\ldots,n)$, since ${\bf r=\bigcup_{i=1}^n r_i}$. So $\mathcal{P}$ is (A)  over $\{\mathcal{P}_i\mid i=1,\ldots,n\}$ if and only if $\mathcal{P}$ is (CA) over $\{\mathcal{P}_i\mid i=1,\ldots,n\}$.

\section{Main results}\label{S:df_in_F}
\begin{thm}\label{th5}
Let ${\bf F}$ be the free group with basis ${\bf x}$, $ \bf{r=\bigcup_{i=1}^n r_i}$ be a set of cyclically reduced words in $\bf{x}\cup\bf{x}^{-1}$, and $\mathcal{P} = \langle \bf{x}\mid \bf{r}\rangle$ be the presentation satisfying the $RH$-hypothesis.
For $i=1,\ldots, n$, let ${\bf R_i}$ be the normal closure of the set ${\bf r_i}$ in ${\bf F}$, $\mathcal{P}_i=\langle{\bf x}\mid{\bf r_i}\rangle$,
${\bf N_i} = \prod_{j\neq i}{\bf R_j}$, $ \bf{\hat{r}_i=\bigcup_{j\neq i} r_j}$, ${\bf R} = \prod_{i=1}^n{\bf R_i}$.

If $\pi_2(\mathcal{P})$ is generated over $\{\mathcal{P}_i\mid i=1,\ldots, n\}$ by a collection $\bf Y$ of based spherical pictures over $\mathcal{P}$, then for $i=1,\ldots, n$, the group $\frac{{\bf R_i}\cap {\bf N_i}}{[{\bf R_i},{\bf N_i}]}$ is generated by
 $$\{[WRW^{-1}][{\bf R_i}, {\bf N_i}]\mid R\in{\bf r_i}\cap{\bf \hat{r}_i}, [W]\in{\bf W}\}\cup\{[WV_YW^{-1}][{\bf R_i}, {\bf N_i}]\mid Y\in{\bf Y}, [W]\in{\bf W}\},$$ where
$V_Y$ is a label of a simple closed path in a based spherical picture $Y$, separating the vertices with $\bf r_i$-labels and the vertices with $(\bf r - r_i)$-labels, ${\bf W}\subseteq {\bf F}$ is a set of representatives of all the cosets of $\bf R$ in ${\bf F}$.
\end{thm}

 It is known \cite{ratcl, DuElGi92, huebsch} that if $\{{\bf R_1},\ldots,{\bf R_n}\}$ is independent, then $\mathcal{P}$ is (A)  over $\{\mathcal{P}_i\mid i=1,\ldots,n\}$. The converse statement for $n=2$ is shown in \cite{ratcl}. From Theorem \ref{th5}, we obtain the following generalization (take $\bf Y$ empty).
\begin{cor}\label{corOFth5}
Let ${\bf F}$ be the free group with basis ${\bf x}$, $ \bf{r=\bigcup_{i=1}^n r_i}$ be a set of cyclically reduced words in $\bf{x}\cup\bf{x}^{-1}$,  where $\bf r_1, \ldots, r_n$ are mutually disjoint sets, and $\mathcal{P} = \langle \bf{x}\mid \bf{r}\rangle$ be the presentation satisfying the $RH$-hypothesis.
For $i=1,\ldots, n$, let ${\bf R_i}$ be the normal closure of ${\bf r_i}$ in ${\bf F}$, $\mathcal{P}_i=\langle{\bf x}\mid{\bf r_i}\rangle$.

If $\mathcal{P}$ is (A) over $\{\mathcal{P}_i\mid i=1,\ldots,n\}$, then $\{{\bf R_1},\ldots,{\bf R_n}\}$ is independent.
\end{cor}

Thus, for $\bf{r=\bigsqcup_{i=1}^n r_i}$, $\{{\bf R_1},\ldots,{\bf R_n}\}$ is independent if and only if $\mathcal{P}$ is (A) over $\{\mathcal{P}_i\mid i=1,\ldots,n\}$.

\section{The proof of  Theorem \ref{th5}}
Let us consider the case $i=1$ (the proof for $i=2,\ldots,n$ is similar).

By $\mathfrak{N}$ denote the normal closure of $\{[R]\mid R\in{\bf r_1}\cap{\bf \hat{r}_1}\}\cup\{ [V_Y]\mid Y\in{\bf Y}\}$ in ${\bf F}$. We need to prove that if $\pi_2(\mathcal{P})$ is generated over $\{\mathcal{P}_j\mid j=1,\ldots, n\}$ by $\bf Y$, then modulo $\mathfrak{N}[{\bf R_1},{\bf N_1}]$, an arbitrary element $[U]\in {\bf R_1}\cap {\bf N_1}$ is equal to the identity.
This will imply that $\frac{{\bf R_1}\cap {\bf N_1}}{[{\bf R_1},{\bf N_1}]}$ is generated by
$$\{[WRW^{-1}][{\bf R_1}, {\bf N_1}]\mid R\in{\bf r_1}\cap{\bf \hat{r}_1}, [W]\in{\bf W}\}\cup\{[WV_YW^{-1}][{\bf R_1}, {\bf N_1}]\mid Y\in{\bf Y}, [W]\in{\bf W}\},$$ since $[V_Y]\in {\bf R_1}\cap {\bf N_1}$, ${\bf r_1}\cap{\bf \hat{r}_1}\subset {\bf R_1}\cap {\bf N_1}$ and ${\bf R}={\bf R_1}{\bf N_1}$.

Since $[U]$ belongs both to  ${\bf R_1}$ and to ${\bf N_1}$, by Lemma \ref{lemVKampen} there are a picture $\mathbf{P}_{{\bf R_1}}$ over $\mathcal{P}_1=\langle{\bf x}\mid{\bf r_1}\rangle$ with boundary label identically equal to $U$ and a picture $\mathbf{P}_{{\bf N_1}}$ over $\mathcal{P}_{({\bf r - r_1})}=\langle{\bf x}\mid({\bf r - r_1})\rangle$ with boundary label identically equal to $U^{-1}$.
Suppose that $U$ is identically equal to $x_1x_2\ldots x_m$, where $x_j\in {\bf x}\cup {\bf x}^{-1}$.  Then, the arcs $\alpha_1,\ldots,\alpha_m$ met $\partial\mathbf{P}_{{\bf R_1}}$ have the labels respectively $x_1, \ldots,x_m$, and the arcs $\beta_1,\ldots, \beta_m$ met $\partial\mathbf{P}_{{\bf N_1}}$ have the labels respectively $x_m,\ldots, x_1$. Paste $\mathbf{P}_{{\bf R_1}}$ and $\mathbf{P}_{\bf N_1}$ by their boundaries so that for $j=1,\ldots, m$, the arc $\alpha_j$ extends the arc $\beta_{m-(j-1)}$ and the global basepoint $O_{{\bf R_1}}$ of $\mathbf{P}_{{\bf R_1}}$ coincides with the global basepoint $O_{{\bf N_1}}$ of $\mathbf{P}_{\bf N_1}$. In the obtained two-sphere, choose a small closed disk $D$ in the interior of any region of $\mathbf{P}_{{\bf R_1}}$ and cut $D-{\partial D}$ out it to get a spherical picture  $\mathbf{P}$ over $\mathcal{P} = \langle \bf{x}\mid \bf{r}\rangle$ on the disk $D^2$ with boundary $\partial\mathbf{P} (= \partial D^2)= \partial D$.

The pasted boundaries $\partial\mathbf{P}_{{\bf R_1}}$ and $\partial\mathbf{P}_{{\bf N_1}}$ give a path $Equ$ on $D^2$. We will call $Equ$ the {\it equator}. The pasted points $O_{{\bf R_1}}$ and $O_{{\bf N_1}}$ give a point $p\in Equ$.
Fix an orientation of $Equ$ so that the label of $Equ-\{p\}$ is equal identically to $U$. Below
the label of $Equ-\{p\}$ under this orientation will be called the {\it equatorial label}.
The part of $\mathbf{P}$ corresponding to $\mathbf{P}_{{\bf R_1}}$ (resp., $\mathbf{P}_{{\bf N_1}}$) will be called the {\it $\bf r_1$-hemisphere} (resp., the {\it $({\bf r-r_1})$-hemisphere}).

 Below we will transform $\mathbf{P}$ by planar isotopy, $BRIDGE$, $FLOAT^{\pm 1}$, $FOLD^{\pm 1}$, $DELETE({\bf X})^{\pm 1}$ under the following conditions: the equatorial label is not changed modulo $\mathfrak{N}[{\bf R_1},{\bf N_1}]$, all vertices of $\mathbf{P}$ with $\bf r_1$-labels remain only in the $\bf r_1$-hemisphere, all vertices with $({\bf r-r_1})$-labels remain only in the $({\bf r-r_1})$-hemisphere. Deformations (operations), satisfied these conditions, will  be called {\it admissible}. A picture, in which the equator divides the vertices under these conditions, will be called a {\it picture with equator}.
The aim of these operations is to reduce the picture $\mathbf{P}$ with equator to a picture with equator with boundary label equal to the identity in the free group, that will imply that the initial equatorial label, i.e. $[U]$, belongs to $\mathfrak{N}[{\bf R_1},{\bf N_1}]$.
To proof the existence of the desired operations, we will need two auxiliary statements: Lemma \ref{lem_help1} and Lemma \ref{lem_help2}.
 In these lemmas, we use the notation of Theorem \ref{th5}, we let  $\bf X$ denote the union ${\bf Y}\cup\bigcup_{i= 1}^n{\bf X}_i$, where ${\bf X}_i$ is the collection of all based spherical picture over $\mathcal{P}_i$, $i=1,\ldots,n$; by an {\it ${\bf Z}^W$-picture} we mean a spherical picture, containing only one $\bf Z$-picture and, possibly, closed arcs encircling this $\bf Z$-picture, where $\bf Z$ is a given collection of spherical pictures.

\begin{lem} \label{lem_help1}
If $\pi_2(\mathcal{P})$ is generated over $\{\mathcal{P}_j\mid j=1,\ldots, n\}$ by $\bf Y$,
then the picture $\mathbf{P}$ with equator $Equ$ can be reduced by a finite number of admissible operations to a picture $\widetilde{\mathbf{P}}$ with equator, being a finite sum of ${\bf X}^W$-pictures and, possibly, also containing some closed arcs encircling the point $p$.
\end{lem}
{\bf Proof of Lemma \ref{lem_help1}.}

Since $\pi_2(\mathcal{P})$ is generated over $\{\mathcal{P}_j\mid j=1,\ldots, n\}$ by $\bf Y$, then by Theorem \ref{thmPr3}, the based spherical picture $\mathbf{P}$ over $\mathcal{P}$ is $\bf X$-equivalent to the empty picture, i.e., there are a finite sequence of based spherical pictures $\mathbf{P}_0, \mathbf{P}_1,\ldots, \mathbf{P}_s$ and a finite sequence $\mathfrak{f}_1,\ldots, \mathfrak{f}_s$ of operations $BRIDGE$, $FLOAT^{\pm 1}$, $FOLD^{\pm 1}$, $DELETE({\bf X})^{\pm 1}$, which transforms $\mathbf{P}$ (up to planar isotopy) to the empty picture so that $\mathfrak{f}_j: \mathbf{P}_{j-1} \mapsto \mathbf{P}_j$, $j=1,\ldots,s$, $\mathbf{P}_0 = \mathbf{P}$, $\mathbf{P}_s$ is the empty picture.

Note that any folding pair is a spherical picture over some presentation $\mathcal{P}_i$, $i=1,\ldots,n$, and, hence, it belongs to ${\bf X}_i$.
In this case $FOLD^{\pm 1}$ is a special case of $DELETE(\bf{X})^{\pm 1}$, therefore, below $FOLD^{\pm 1}$ will not be considered separately.

Since $\mathfrak{f}_1,\ldots, \mathfrak{f}_s$ are not necessarily admissible, using the sequences $\mathbf{P}_0, \mathbf{P}_1,\ldots, \mathbf{P}_s$ and $\mathfrak{f}_1,\ldots, \mathfrak{f}_s$, we will construct two new sequences  $\mathfrak{Z}_0, \mathfrak{Z}_1,\ldots,\mathfrak{Z}_s$ and $\mathfrak{g}_1,\ldots, \mathfrak{g}_s, \mathfrak{g}_{s+1}$, where, for $i=1,\ldots,s$, $\mathfrak{Z}_i$ is a collection of disjoint spherical subpictures in $\mathbf{P}_i$, not containing the whole of the equator,  and $\mathfrak{g}_i$ is an admissible operation. In addition, the sequence $\mathfrak{g}_1,\ldots, \mathfrak{g}_s, \mathfrak{g}_{s+1}$ will transform $\mathbf{P}$, as a picture with equator, to $\widetilde{\mathbf{P}}$ so that  $\mathfrak{g}_j: \widetilde{\mathbf{P}}_{j-1} \mapsto \widetilde{\mathbf{P}}_j$, where $\widetilde{\mathbf{P}}_0 = \mathbf{P}$, $\widetilde{\mathbf{P}}_j = \mathbf{P}_j\cup\mathfrak{Z}_j$, $\widetilde{\mathbf{P}}_{s+1} =\widetilde{\mathbf{P}}$ are pictures with equator, $j=1,\ldots,s$. This will prove Lemma \ref{lem_help1}.

We let $\mathfrak{Z}_j^0$, $j=1,\ldots,s$, denote a finite collection of disjoint disks in the interior of $D^2$, containing the empty spherical subpictures obtained from subpictures of $\mathfrak{Z}_j$ by deletion of all arcs and all vertices.

In the case of a planar isotopy, we may assume that this isotopy reduces $\mathbf{P}_{j-1}$ to
$\mathbf{P}_{j}=F_1(\mathbf{P}_{j-1})$, where $F_t : D^2\times [0,1] \to
D^2\times [0,1]$ is a continuous isotopy of the disk $D^2$, so that
\begin{itemize}
\item[(i)] $F_t$ leaves all vertices fixed, i.e., for any $t\in [0,1]$ and each vertex $\Delta$,
$F_t(\Delta)=\Delta$;
\item[(ii)] for any $t\in [0,1]$ and each arc $\alpha$, the intersection of the arc $F_t(\alpha)$ and the equator $Equ$ consists of a finite number of points; moreover, if $Equ$ intersects
$F_1(\alpha)$, then it crosses it, and doesn't just touch it;
\item[(iii)] for any arc $\alpha$,  the arc $F_1(\alpha)$ does not intersect any disk of $\mathfrak{Z}_{j-1}^0$.
\end{itemize}

 If for any $t\in (0,1)$ and each arc $\alpha$, the arc $F_t(\alpha)$ does not intersect any disk of $\mathfrak{Z}_{j-1}^0$, then this isotopy is called  {\it admissible}, otherwise it is called {\it inadmissible}. An admissible isotopy does not change the equatorial label, as an element of the free group. Since an admissible isotopy is an admissible operation, below we will assume operations to be equal if they are equal up to admissible isotopy. Include operations of inadmissible isotopy in the list of operations $\mathfrak{f}_1,\ldots, \mathfrak{f}_s$.

 We will construct $\mathfrak{g}_j$, $j=1,\ldots,s$, so that  $\mathfrak{g}_j$ transforms
 $(\mathbf{P}_{j-1}-\mathfrak{Z}_{j})$ to $(\mathbf{P}_j-\mathfrak{Z}_j)$ in just the same way (up to admissible isotopy of $(\mathbf{P}_j-\mathfrak{Z}_j)$) as $\mathfrak{f}_j$ transforms $(\mathbf{P}_{j-1}-\mathfrak{Z}_{j}^0)$ to $(\mathbf{P}_j-\mathfrak{Z}_j^0)$. In addition we can always assume that the arcs of $\widetilde{\mathbf{P}}_j$ and the boundary of any subpicture from $\mathfrak{Z}_j^0$ intersect the equator not more than finitely many times.

As $\mathfrak{Z}_0$, take a set of a single spherical subpicture in $\mathbf{P} = \mathbf{P}_0$ such that this subpicture is empty and contains the point $p\in Equ$  and a connected part of the equator.

Let us define $\mathfrak{Z}_j$ and $\mathfrak{g}_j$, $j=1,\ldots,s$, by induction on $j$.

Assume that $\mathfrak{Z}_0, \mathfrak{Z}_1,\ldots,\mathfrak{Z}_{j-1}$ and $\mathfrak{g}_1,\ldots, \mathfrak{g}_{j-1}$ have been already defined. Construct $\mathfrak{Z}_{j}$ and $\mathfrak{g}_{j}$ by $\mathfrak{f}_j: \mathbf{P}_{j-1} \mapsto \mathbf{P}_j$ as follows.  The operation $\mathfrak{f}_j$ is one of the operations: an inadmissible isotopy, $BRIDGE$, $FLOAT^{\pm 1}$, $DELETE({\bf X})^{\pm 1}$.

 \underline{Case 1. Inadmissible isotopy.}

 There are an arc $\alpha$ (labeled by $x\in \bf{x}\cup\bf{x}^{-1}$) in $\mathbf{P}_{j-1}$ and $t\in (0,1)$ so that $F_t(\alpha)$ intersects some empty spherical subpicture $\mathbf{P}_{\xi}^0$ from $\mathfrak{Z}_{j-1}^0$. Denote by $\mathbf{P}_{\xi}$ the spherical subpicture from $\mathfrak{Z}_{j-1}$, which $\mathbf{P}_{\xi}^0$ is obtained from. To obtain $\mathfrak{Z}_j$ from $\mathfrak{Z}_{j-1}$, (for each arc $\alpha$, each such subpicture $\mathbf{P}_{\xi}^0$ and each passing of $F_t(\alpha)$ through $\mathbf{P}_{\xi}^0$) add to $\mathbf{P}_{\xi}$ a closed arc (labeled by $x$), encircling all objects (that may be arcs, vertices, the point $p$) being already in $\mathbf{P}_{\xi}$ . The action of $\mathfrak{g}_j$ on
 $(\mathbf{P}_{j-1}-\mathfrak{Z}_{j})$ coincides with the action of $\mathfrak{f}_j$ on $(\mathbf{P}_{j-1}-\mathfrak{Z}_{j}^0)$. The action of $\mathfrak{g}_j$ on $\mathfrak{Z}_{j}$ corresponds to the operation $BRIDGE$ performed on the arc $\alpha$ in the interior of the disk of $\mathbf{P}_{\xi}$.
 So $\mathfrak{g}_j$ does not change the equatorial label as an element of the free group. See Figure 3.

\unitlength=1mm
\unitlength=1mm
\begin{picture}(120,50)(-5,0)

\qbezier[30](42,33)(50,40)(58,33)
\qbezier[30](42,17)(50,10)(58,17)
\qbezier[30](42,33)(35,25)(42,17)
\qbezier[30](58,33)(65,25)(58,17)

       \put(48,10){$\mathbf{P}_{\xi}^0$}

\qbezier[30](92,33)(100,40)(108,33)
\qbezier[30](92,17)(100,10)(108,17)
\qbezier[30](92,33)(85,25)(92,17)
\qbezier[30](108,33)(115,25)(108,17)

       \put(98,10){$\mathbf{P}_{\xi}^0$}

       \thicklines

       \put(70,25){\vector(1,0){10}}
       \put(73,28){\rm $\mathfrak{f}_j$}

       \thinlines
       \put(30,40){\line(0,-1){30}}

       \put(28,25){\vector(1,0){4}}

       \put(27,26){\rm $x$}

        \put(120,40){\line(0,-1){30}}
        \put(118,25){\vector(1,0){4}}

       \put(117,26){\rm $x$}
\end{picture}

\unitlength=1mm
\unitlength=1mm
\begin{picture}(120,50)(-5,0)

\qbezier[30](42,33)(50,40)(58,33)
\qbezier[30](42,17)(50,10)(58,17)
\qbezier[30](42,33)(35,25)(42,17)
\qbezier[30](58,33)(65,25)(58,17)

       \put(48,10){$\mathbf{P}_{\xi}$}

\qbezier[30](92,33)(100,40)(108,33)
\qbezier[30](92,17)(100,10)(108,17)
\qbezier[30](92,33)(85,25)(92,17)
\qbezier[30](108,33)(115,25)(108,17)

       \put(98,10){$\mathbf{P}_{\xi}$}

       \put(51,21){\circle*{2}}
        \put(51,29){\circle*{2}}
         \put(47,25){\circle*{2}}
         \thinlines
       \put(47,25){\line(1,1){4}}
       \put(47,25){\line(1,-1){4}}
       \put(51,21){\line(0,1){8}}

       \put(100,25){\circle{20}}

        \put(101,21){\circle*{2}}
        \put(101,29){\circle*{2}}
         \put(97,25){\circle*{2}}
         \thinlines
       \put(97,25){\line(1,1){4}}
       \put(97,25){\line(1,-1){4}}
       \put(101,21){\line(0,1){8}}

       \thicklines

       \put(70,25){\vector(1,0){10}}
       \put(73,28){\rm $\mathfrak{g}_j$}

       \thinlines
       \put(30,40){\line(0,-1){30}}

       \put(28,25){\vector(1,0){4}}

       \put(27,26){\rm $x$}

       \put(109,25){\vector(-1,0){4}}
       \put(108,23){\tiny \rm $x$}

        \put(120,40){\line(0,-1){30}}
        \put(118,25){\vector(1,0){4}}

       \put(117,26){\rm $x$}

\end{picture}

\begin{center} Figure 3
\end{center}

 \underline{Case 2. $FLOAT$, $DELETE({\bf X})$.}

  In $\mathbf{P}_{j-1}$, there is a spherical subpicture $\mathbf{P}_{\eta}$ containing only a floating circle  (resp., only an ${\bf X}$-picture). The operation  $\mathfrak{f}_j$ deletes $\mathbf{P}_{\eta}$ from $\mathbf{P}_{j-1}$.

 Applying an admissible isotopy to $\partial\mathbf{P}_{\eta}$ and $\mathfrak{f}_j$ if necessary, we may assume that $\mathbf{P}_{\eta}$ does not contain the whole of the equator, and $\partial\mathbf{P}_{\eta}$ does not intersect any disk from $\mathfrak{Z}_{j-1}^0$. Note that the arcs and the vertices of  $\mathbf{P}_{\eta}$ do not intersect $\mathfrak{Z}_{j-1}^0$ as well.

 Let $\mathbf{P}_{\xi_1}^0$, ..., $\mathbf{P}_{\xi_m}^0$ be all spherical subpictures from $\mathfrak{Z}_{j-1}^0$ in the interior of the disk of $\mathbf{P}_{\eta}$, and let $\mathbf{P}_{\xi_1}$, ..., $\mathbf{P}_{\xi_m}$ be the corresponding subpictures from $\mathfrak{Z}_{j-1}$.
  The collection $\mathfrak{Z}_{j}$ is obtained from $\mathfrak{Z}_{j-1}$ by deleting $\mathbf{P}_{\xi_1}$, ..., $\mathbf{P}_{\xi_m}$ and adding the spherical subpicture from $\widetilde{\mathbf{P}}_{j-1}$ encircled by the path $\partial\mathbf{P}_{\eta}$ (this subpicture contains the disjoint union of $\mathbf{P}_{\eta}$ and $\mathbf{P}_{\xi_1}$, ..., $\mathbf{P}_{\xi_m}$). The operation $\mathfrak{g}_j$ acts identically (up to admissible isotopy).

\underline{Case 3. $FLOAT^{-1}$, $DELETE({\bf X})^{-1}$.}

Let $\mathfrak{f}_j$ be an operation $FLOAT^{-1}$ (resp., $DELETE({\bf X})^{-1}$).
The operation $\mathfrak{f}_j$ inserts a spherical subpicture $\mathbf{P}_{\eta}$ in $(\mathbf{P}_{j-1} - \mathfrak{Z}_{j-1}^0)$ such that $\mathbf{P}_{\eta}$ contains only a floating circle (resp., only an ${\bf X}$-picture).

Applying an admissible isotopy if necessary, we can obtain the following. If $\mathfrak{f}_j$ is $FLOAT^{-1}$, the subpicture $\mathbf{P}_{\eta}$ does not intersect the equator. If $\mathfrak{f}_j$ is $DELETE({\bf X})^{-1}$ and $\mathbf{P}_{\eta}$ contains the vertices only with $\bf r_1$-labels or only with $(\bf r - r_1)$-labels, then the whole of $\mathbf{P}_{\eta}$ is in the hemisphere corresponding to the labels of the vertices in $\mathbf{P}_{\eta}$. If $\mathfrak{f}_j$ is $DELETE({\bf X})^{-1}$ and $\mathbf{P}_{\eta}$ contains the vertices both with $\bf r_1$-labels and with $(\bf r - r_1)$-labels (we may assume that $\mathbf{P}_{\eta}\in{\bf Y}$), then $\mathbf{P}_{\eta}$ is disposed so that  $\hat{Equ} = \mathbf{P}_{\eta}\cap Equ$ is connected and the equator divides the vertices of $\mathbf{P}_{\eta}$ into two parts: the vertices with $\bf r_1$-labels (in the $\bf r_1$-hemisphere) and the vertices with $(\bf r - r_1)$-labels (in the $(\bf r - r_1)$-hemisphere), in addition the label $U_{\eta}$ of the path $\hat{Equ}$ is such that $[U_{\eta}]\in \mathfrak{N}$.

 Put $\mathfrak{Z}_{j} =\mathfrak{Z}_{j-1}$. The action of $\mathfrak{g}_j$ on
 $(\mathbf{P}_{j-1}-\mathfrak{Z}_{j})$ coincides with the action of $\mathfrak{f}_j$ on $(\mathbf{P}_{j-1}-\mathfrak{Z}_{j}^0)$, and the action of $\mathfrak{g}_j$ on $\mathfrak{Z}_{j}$ is identical.
 The operation $\mathfrak{g}_j$ does not change the equatorial label modulo $\mathfrak{N}[{\bf R_1},{\bf N_1}]$.

\underline{Case 4. $BRIDGE$.}

We may assume that $\mathfrak{f}_j$ acts only on $(\mathbf{P}_{j-1} - \mathfrak{Z}_{j-1}^0)$.

Put $\mathfrak{Z}_{j} =\mathfrak{Z}_{j-1}$. The action of $\mathfrak{g}_j$ on
 $(\mathbf{P}_{j-1}-\mathfrak{Z}_{j})$ coincides with the action of $\mathfrak{f}_j$ on $(\mathbf{P}_{j-1}-\mathfrak{Z}_{j}^0)$, and the action of $\mathfrak{g}_j$ on $\mathfrak{Z}_{j}$ is identical. The operation $\mathfrak{g}_j$ does not change the equatorial label as an element in the free group.

 To complete the proof of Lemma \ref{lem_help1}, it remains to define  $\mathfrak{g}_{s+1}: \widetilde{\mathbf{P}}_{s} \mapsto \widetilde{\mathbf{P}}_{s+1}$, where  $\widetilde{\mathbf{P}}_s = \mathbf{P}_s\cup\mathfrak{Z}_s$, $\widetilde{\mathbf{P}}_{s+1} =\widetilde{\mathbf{P}}$.
 Since $\mathbf{P}_s$ is the empty picture, the arcs and the vertices of $\widetilde{\mathbf{P}}_s$ belong to spherical subpictures from $\mathfrak{Z}_{s}$.
 The operation $\mathfrak{g}_{s+1}$ transforms each spherical picture  $\mathbf{P}_{\mu}$ from $\mathfrak{Z}_{s}$ as follows.

By construction of $\mathfrak{Z}_{s}$,
 $\mathbf{P}_{\mu}$ is the union of embedded one in other spherical pictures each of which either is an ${\bf X}^W$-picture, or contains only closed arcs and, possibly, the point $p$.
The operation $\mathfrak{g}_{s+1}$ decomposes $\mathbf{P}_{\mu}$ as a sum of ${\bf X}^W$-pictures and spherical pictures, containing only closed arcs and, possibly, $p$, by applying admissibly isotopy and $BRIDGE$ (see an example on Figure 4) so that the vertices remain in their own hemispheres, no summand contains the whole of the equator, and the boundary of each summand  intersects the equator not more than finitely many times.
After that $\mathfrak{g}_{s+1}$ deletes all floating circles not encircling $p$ by applying $FLOAT$.

The operation $\mathfrak{g}_{s+1}$ does not change the equatorial label as an element of the free group.

\unitlength=1mm
\begin{picture}(125,50)(1,0)
       \thicklines
       \put(77,25){\vector(1,0){10}}
       \thinlines
       \multiput(10,30)(105,0){2}{\circle*{2}}
       \multiput(50,30)(105,0){2}{\circle*{2}}
       \multiput(10,40)(105,0){2}{\circle*{2}}
       \multiput(50,40)(105,0){2}{\circle*{2}}

       \multiput(10,31)(105,0){2}{\line(0,1){8}}
       \multiput(50,31)(105,0){2}{\line(0,1){8}}
       \multiput(11,40)(105,0){2}{\line(1,0){38}}
       \put(30,29){\oval(40,38)[b]}

       \multiput(30,22)(105,-2){2}{\circle{2}}
       \multiput(26,15)(105,-2){2}{\circle{2}}
       \multiput(34,15)(105,-2){2}{\circle{2}}

       \multiput(31,21)(105,-2){2}{\line(1,-2){2.5}}
       \multiput(29,21)(105,-2){2}{\line(-1,-2){2.5}}
       \multiput(27,15)(105,-2){2}{\line(1,0){6}}

       \put(115,30){\line(1,0){40}}
       \put(145,28){\vector(0,1){5}}

       \put(135,16.5){\oval(40,13)}

      \put(40,8){\vector(0,1){5}}
       \put(145,25){\vector(0,-1){5}}

       \put(146,20){\rm x}
       \multiput(41,11)(105,20){2}{\rm x}
       \multiput(0,26)(107,8){2}{\rm $\mathbf{P}'$}
       \multiput(22,18)(105,-2){2}{\rm $\mathbf{P}''$}

\qbezier[60](8,42)(30,55)(52,42)
\qbezier[60](8,10)(30,-3)(52,10)
\qbezier[60](8,42)(-11,26)(8,10)
\qbezier[60](52,42)(71,26)(52,10)

\qbezier[60](113,42)(135,48)(157,42)
\qbezier[60](113,28)(135,22)(157,28)
\qbezier[40](113,42)(98,35)(113,28)
\qbezier[40](157,42)(172,35)(157,28)

\qbezier[30](23,22)(30,26)(37,22)
\qbezier[30](23,14)(30,10)(37,14)
\qbezier[20](23,22)(18,18)(23,14)
\qbezier[20](37,22)(42,18)(37,14)

\qbezier[30](128,20)(135,24)(142,20)
\qbezier[30](128,12)(135,9)(142,12)
\qbezier[20](128,20)(123,16)(128,12)
\qbezier[20](142,20)(147,16)(142,12)
\end{picture}

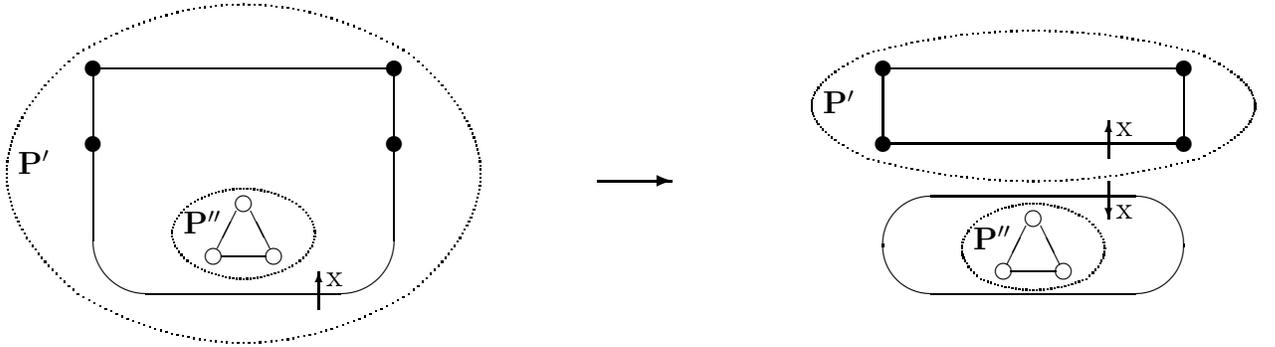
\begin{center} Figure 4. {Partition of subpictures $\mathbf{P}'$ and $\mathbf{P}''$}\\{\tiny (Not to complicate the figure, the orientation and the label are indicated only for the arc which is transformed by $BRIDGE$.)}
\end{center}

 \edvo

\begin{lem} \label{lem_help2}

If a based spherical picture $\tilde{\mathbf{P}}$ with the equator $Equ$ is a finite sum of ${\bf X}^W$-pictures and, possibly, contains some closed arcs around the point $p$, then by a finite number of admissible operations, $\tilde{\mathbf{P}}$
can be reduced to a picture with equator with the equatorial label equal to the identity in the free group.
\end{lem}
{\bf Proof of Lemma \ref{lem_help2}.}
Below we will use the operation $COMMUTE$ depicted on Figure 5. This operation changes the equatorial label on an element from $[{\bf R_1},{\bf N_1}]$ (see details in \cite{ok1}). It can be realized as a planar isotopy and a finite number of $DELETE({\bf X})^{\pm 1}$, $BRIDGE$, $FLOAT^{\pm 1}$.

\unitlength=1mm
\begin{picture}(120,50)(-9,0)
       \thicklines
       \multiput(30,10)(10,0){3}{\line(1,0){5}}
       \multiput(85,10)(10,0){3}{\line(1,0){5}}
       \multiput(35,40)(10,0){3}{\line(1,0){5}}
       \multiput(90,40)(10,0){3}{\line(1,0){5}}

       \multiput(30,15)(0,10){3}{\line(0,1){5}}
       \multiput(85,15)(0,10){3}{\line(0,1){5}}
       \multiput(60,10)(0,10){3}{\line(0,1){5}}
       \multiput(115,10)(0,10){3}{\line(0,1){5}}

       \put(67,25){\vector(1,0){10}}

        \thicklines
        \multiput(30,25)(55,0){2}{\line(1,0){30}}
        \thinlines
        \multiput(20,25)(97,0){2}{${Equ}$}

        \put(38,10){\line(0,1){20}}
        \put(38,31){\circle{2}}
        \qbezier[10](35,30)(35,34)(38,34)
        \qbezier[10](41,30)(41,34)(38,34)
        \qbezier[30](35,10)(35,20)(35,30)
        \qbezier[30](41,10)(41,20)(41,30)

        \put(52,40){\line(0,-1){20}}
        \put(52,19){\circle*{2}}
        \qbezier[10](49,20)(49,16)(52,16)
        \qbezier[10](55,20)(55,16)(52,16)
        \qbezier[30](49,40)(49,30)(49,20)
        \qbezier[30](55,40)(55,30)(55,20)

        \put(93,25){\line(0,-1){5}}
        \put(93,19){\circle*{2}}
        \qbezier[10](90,20)(90,16)(93,16)
        \qbezier[10](96,20)(96,16)(93,16)
        \qbezier[10](90,20)(90,22)(90,25)
        \qbezier[10](96,20)(96,22)(96,25)

        \put(93,10){\line(0,1){2}}
        \qbezier[30](94,15)(105,17)(104,25)
        \qbezier[10](90,10)(89,14)(94,15)
        \qbezier(96,14)(107,14)(107,25)
        \qbezier(93,10)(92,13)(96,14)
        \qbezier[30](99,12)(110,13)(110,25)
        \qbezier[5](96,10)(95,12)(99,12)

        \put(107,25){\line(0,1){5}}
        \put(107,31){\circle{2}}
        \qbezier[10](104,30)(104,34)(107,34)
        \qbezier[10](110,30)(110,34)(107,34)
        \qbezier[10](104,25)(104,28)(104,30)
        \qbezier[10](110,25)(110,28)(110,30)

        \put(107,40){\line(0,-1){2}}

        \qbezier[30](106,35)(95,33)(96,25)
        \qbezier[10](110,40)(111,36)(106,35)
        \qbezier(104,36)(93,36)(93,25)
        \qbezier(107,40)(108,37)(104,36)
        \qbezier[30](101,38)(90,37)(90,25)
        \qbezier[5](104,40)(105,38)(101,38)

\end{picture}

\begin{center}
Figure 5
\end{center}

At first, applying $FLOAT$ and $DELETE({\bf X})$ to $\widetilde{\mathbf{P}}$, we delete all ${\bf X}^W$-pictures, not intersecting $Equ$. This does not change the equatorial label.
Further, by means of planar isotopy and $COMMUTE$, we obtain that for each ${\bf X}^W$-picture $\mathbf{P}_{\eta}$, the intersection $Equ\cap \mathbf{P}_{\eta}$ is connected, i.e., $Equ$ divides $\mathbf{P}_{\eta}$ into two parts: a subpicture over $\mathcal{P}_1=\langle{\bf x}\mid{\bf r_1}\rangle$ in the $\bf r_1$-hemisphere and a subpicture over $\mathcal{P}_{({\bf r - r_1})}=\langle{\bf x}\mid({\bf r - r_1})\rangle$ in the $(\bf r - r_1)$-hemisphere.
If at least one of these parts does not contain vertices, the label of $Equ\cap \mathbf{P}_{\eta}$ is equal to the identity in the free group.
Otherwise either $\mathbf{P}_{\eta}$ is an ${\bf Y}^W$-picture, or $\mathbf{P}_{\eta}$ contains vertices with labels from ${\bf r_1}\cap{\bf \hat{r}_1}$. In the first case the label of $Equ\cap \mathbf{P}_{\eta}$ is equal to $[WV_YW^{-1}]$ in the free group, where
$V_Y$ is the label of a simple closed path in a based spherical picture $Y\in {\bf Y}$, separating the vertices with $\bf r_1$-labels and the vertices with $(\bf r - r_1)$-labels. In the second case the label of $Equ\cap \mathbf{P}_{\eta}$ is equal to the product of elements of the form
 $[WR^{\pm 1}W^{-1}]$, where $R\in{\bf r_1}\cap{\bf \hat{r}_1}$. In each of these cases the label of $Equ\cap \mathbf{P}_{\eta}$ belongs to $\mathfrak{N}[{\bf R_1},{\bf N_1}]$.
 Now we delete all ${\bf X}^W$-pictures from $\widetilde{\mathbf{P}}$. It follows from the above arguments that this operation is admissible.
After such operation only several closed arcs encircling $p$ may remain in $\widetilde{\mathbf{P}}$. So the equatorial label is equal to the identity in the free group. This completes the proof of Lemma \ref{lem_help2}. \edvo

Let us continue the proof of Theorem \ref{th5}. By Lemma \ref{lem_help1}, the picture $\mathbf{P}$ with the equator can be reduced  by a finite number of admissible operations to a picture $\widetilde{\mathbf{P}}$ with equator, being a finite sum of ${\bf X}^W$-pictures and, possibly, also containing some closed arcs around $p$.
By Lemma \ref{lem_help2}, $\tilde{\mathbf{P}}$ can be reduced by a finite number of admissible operations
to a picture with equator with the equatorial label equal to the identity in the free group. So the equatorial label $[U]$ of the initial picture $\mathbf{P}$ is equal to the identity modulo $\mathfrak{N}[{\bf R_1},{\bf N_1}]$. \edvo

{\it Keywords:} Presentations, subpresentations, asphericity, independent families of normal subgroups, intersection of normal subgroups, mutual commutant
of normal subgroups, spherical pictures, identity sequences.
\newpage

\end{document}